\documentclass[a4paper,12pt]{article}
 \usepackage{amssymb,latexsym,amsmath,amsthm}
 \theoremstyle{definition}
\newtheorem{opr}{Definition}
 \theoremstyle{plain}
\newtheorem{thm}{Theorem}
\newtheorem{lem}{Lemma}
\newtheorem{sled}{Corollary}
\theoremstyle{remark}
 \newtheorem{note}{Remark}

\DeclareMathOperator{\Ker}{Ker}

\DeclareMathOperator{\ad}{ad}

\newcommand{\bl}{\square}
\newcommand{\A}{{\mathcal{A}}(n)}
\newcommand{\C}{\mathbb{C}}
\newcommand{\R}{\mathbb{R}}
\newcommand{\doc}{\emph{Proof. }}
\renewcommand{\leq}{\leqslant}
\renewcommand{\geq}{\geqslant}

\newcommand{\ph}{\varphi}
\newcommand{\itd}[2]{#1_1,\dotsc, #1_{#2}}
\newcommand{\matr}[4]{ \begin{pmatrix} #1 & #2 \\ #3 & #4   \end{pmatrix} }
\newcommand{\gf}{{\mathfrak g}_\varPhi}

\begin{document}
 
  \title {On the dimensions of commutative subalgebras and subgroups}
\date{}
\author{M.V. Milentyeva\thanks{This research is partially supported by RFBR (no.~05-01-00895.)}}
\maketitle
 
{\footnotesize We consider the functions that bound the dimensions of finite-dimensional associative or Lie algebras  in terms of  the dimensions of  their commutative subalgebras.
It is proved that these functions  have quadratic growth. As a result,  we also get similar estimates for the dimension of a Lie group with bounded dimensions of its abelian Lie subgroups.
}

\pagestyle{plain}

\section*{Introduction}

We begin with   definitions of  functions we shall consider.
 \begin{opr}

Let  $A$ be an associative or Lie algebra. We say that  $A$  satisfies a condition ${\mathcal{A}}(n)$ if the dimensions of all commutative subalgebras of $A$ are at most $n$.
\end{opr}

  \begin{opr}
For any integer $n$  we denote by $l_{K}(n)$  ($a_{K}(n)$)  the greatest integer $h$ 
such that there exists a Lie algebra   (an associative algebra)  of dimension $h$ over a field $K$ such that  this algebra satisfies the condition $\A$.
  \end{opr}
  \begin{opr}
We denote by $a^{1}_{K}(n)$ the minimal function such that $\dim A \leq a^{1}_{K}(n)$ for any unital associative algebra $A$ over a field K whose unital commutative subalgebras all have dimensions at most $n$.
  \end{opr}

  \begin{opr}
Let $K$ be the complex number field or the field of real numbers. For any integer $n$  we denote by $g_{K}(n)$    the greatest integer $h$ satisfying the following property: there exists a Lie group of dimension $h$ over  $K$ such that  the dimensions of all commutative Lie subgroups of this group are at most $n$. 
  \end{opr}
Throughout this paper we assume that all algebras are finite dimensional.
 
Formally, it does not follow from the Definitions that this functions are well defined, since they are not necessary finite. Nevertheless, we shall show that they are finite in the case of the complex or real number field and in some other cases, and, moreover, they have quadratic growth. Namely  the following theorem holds.\\ 

\noindent \textbf{Main Theorem.} \emph{The functions we defined have quadratic growth.
Namely the following inequalities hold:\begin{align*}
\frac{n^2+4n-5}{8}&\leq l_{\C}(n) \leq \frac{n^2+17n}{2};\\ \frac{n^2+4n-5}{8}&\leq g_{\C}(n) \leq \frac{n^2+17n}{2};\\  \frac{n^2+4n-5}{8}&\leq a_{\C}(n)\leq \frac {n^2} 2 +5n;\\ 2n^2+n &\leq l_{\R}(n)\leq 4n^2+18n;\\ 2n^2+n &\leq g_{\R}(n)\leq 4n^2+18n;\\ \frac{n^2+4n-5}{8} &\leq a_{\R}(n)\leq \frac{n^2 }{2}+5n;
\end{align*}
if  $K$ is a field of characteristic $0$, then
$$a_K(n) \leq \frac {3n^2+n} 2,$$
$$a_K^1(n) \leq \frac {3n^2+n} 2;$$
 if a field $K$ is algebraically closed, then
$$a_K(n) \leq \frac {n^2} 2+5n,$$
$$a^1_K(n) \leq \frac {n^2} 2+5n;$$
for any field $K$  }
\begin{align*}
l_K(n) &\geq \frac{n^2+4n-5}{8},  \\ a_K(n) &\geq \frac {n^2+4n-5} 8, \\ a^1_K(n) &\geq \frac {n^2+2n } 8.
\end{align*}

Explicit forms for the functions $l_K,\ a_K, \ a^1_K$, and $g_K $ are not known.
However, the Theorem rises the question about asymptotic behavior of  these functions: find the limits  $\frac{f(n)}{n^2}$ as $n\rightarrow \infty$.
 
This work continues the research started in  A.Yu.~Olshanskii's paper~\cite{O78} and in author's one~\cite{M4}, where similar functions for finite $p$-groups and finitely generated torsion-free nilpotent groups, respectively, were studied. We use certain results of these papers to obtain the lower bounds. In fact, the results from [3] and [4] are already
enough to establish the quadratic estimates even for nilpotent algebras and groups of class $2$. To achieve the upper bounds, we will utilize standard 
properties of Lie groups and algebras.

\section{Quadratic upper bounds}
\subsection{ Complex Lie algebras}

\begin{lem} \label {T:A}
Let ${\mathfrak g}$ be a nilpotent Lie algebra over a field $K$. Suppose that $A$ satisfies the condition $\A$; then
 $$ \dim {\mathfrak g} \leq \frac{n(n+1)}{2}.$$
\end{lem}
\doc
Let ${\mathfrak i}\lhd {\mathfrak g}$ be a maximal abelian ideal. Put $s=\dim {\mathfrak i}$.
Denote by $\rho$ the representation of ${\mathfrak g}$ that takes each $x\in{\mathfrak g}$ to the restriction of adjoint action of $x$ to ${\mathfrak i}$ ($\rho(x)=\ad x|_{\mathfrak i}$).
 Let us show that ${\Ker \rho ={\mathfrak i}}$. Obviously, since ${\mathfrak i} $ is abelian, ${\mathfrak i} \subseteq \Ker \rho$. If ${\mathfrak i} \subsetneq \Ker \rho$, then $\Ker \rho / {\mathfrak i}$ is  a nontrivial ideal of the nilpotent factor algebra ${\mathfrak g}/{\mathfrak i}$. Consequently the intersection of $\Ker \rho / {\mathfrak i}$ with the centre of ${\mathfrak g}/{\mathfrak i}$ is nontrivial.
That is there exists $x \in \Ker \rho  \setminus {\mathfrak i}$ such that $[x,y]\in {\mathfrak i}$ for any $y\in {\mathfrak g}$.
And we see that ideal generated by $x$ and ${\mathfrak i}$ is abelian and it contains ${\mathfrak i}$ as a proper subset contradicting our assumption that ${\mathfrak i}$ is maximal.

Since ${\mathfrak g}$ is nilpotent, for any $x\in {\mathfrak g}$ an endomorphism $\ad x$ is also nilpotent.
It follows that the factor algebra ${\mathfrak g}/{\mathfrak i}$ is isomorphic to a   subalgebra of ${\mathfrak gl}_s(K)$ whose elements all are nilpotent. By Engel's theorem, ${\mathfrak g}/{\mathfrak i}$ is isomorphic to a subalgebra that consists of strictly upper  triangular matrices.
Hence $\dim {\mathfrak g}/{\mathfrak i} \leq \frac{s(s-1)}{2} $, so that $\dim {\mathfrak g} \leq \frac{s(s+1)}{2} .$~$\bl$
 
\begin{lem} \label{T:R}
 If ${\mathfrak g}$ is solvable Lie algebra over an algebraically  closed field $K$ of characteristic 0 such that ${\mathfrak g}$ satisfies the condition $\A$, then
 $$ \dim {\mathfrak g} \leq \frac{n(n+3)}{2}.$$
\end{lem}
\doc
We choose a maximal abelian ideal ${\mathfrak i}\lhd {\mathfrak g}$ and put $s=\dim {\mathfrak i}$ as in the proof of the previous Lemma. For each $x\in{\mathfrak g}$ put $\rho(x)=\ad x|_{\mathfrak i}$.
 
Evidently, ${\mathfrak i} \subseteq \Ker \rho$. If ${\mathfrak i} \neq \Ker \rho$, then, by Corollary of Lie's theorem (see Corollary~3 of Theorem~I.5.1,~\cite{B76}), the nontrivial ideal $\Ker \rho / {\mathfrak i} \lhd {\mathfrak g}/ {\mathfrak i}$ contains a  one-dimensional ideal ${\mathfrak j}/ {\mathfrak i}$. But then $ {\mathfrak j}$ is an abelian ideal of $ {\mathfrak g}$ and it is greater than $ {\mathfrak i}$. We have a contradiction.

Therefore ${\Ker \rho ={\mathfrak i}}$ and ${\mathfrak g}/ {\mathfrak i}$ is isomorphic to a  solvable subalgebra of ${\mathfrak gl}_s(K)$. By Lie's theorem, all elements of this subalgebra can be simultaneously represented by upper triangular matrices. Whence $\dim {\mathfrak g}/{\mathfrak i} \leq \frac{s(s+1)}{2} $, so that  $\dim {\mathfrak g} \leq \frac{s(s+3)}{2} .$~$\bl$

\begin{lem} \label{T:B}
Let ${\mathfrak g}$ be a complex simple Lie algebra satisfying the  condition $\A$. Then
\begin{equation} \label{E:A}
   \dim {\mathfrak g} \leq 7n.
\end{equation}
\end{lem}
\doc 
According to the classification  (see Theorems~4.2.13, 4.3.1 and 4.3.3, \cite{VO95}), there are the following types of complex simple Lie algebras (see the first table row):
\begin{center}
\begin{tabular}{|c|c|c|c|c|c|c|c|c|c|}
 \hline
Type of ${\mathfrak g}$& $E_6$&$E_7$&$E_8$&$F_4$&$G_2$&$A_l\ (l\geq 1)$&$B_l\ (l\geq 3)$&$C_l\ (l\geq 2)$ &$D_l\ (l\geq 4)$\\ \hline

 $\dim {\mathfrak g}$ &78 &133& 248 &52&14
 & $l^2+2l $& $2l^2+l $&$2l^2+l$ & $2l^2-l$ \\ \hline
$\dim {\mathfrak h}$&16&27&36&9&3 &$\left[\frac{(l+1)^2}{4}\right]$&$ \frac{l(l-1)}{2}+1$&$\frac{l(l+1)} 2$ &$\frac{l(l-1)} 2$\\ \hline
\end{tabular}
\end{center}
Here index stands for the rank of the corresponding Lie algebra, that is the dimension of a Cartan subalgebra.
The second table row contains the dimensions of the corresponding Lie algebras in accordance with Table~1,~\cite{VO95}.
A.I.~Malcev determined in  \cite{M45} commutative subalgebras of maximal dimensions for simple complex  Lie algebras.
According to this paper, ${\mathfrak g}$ contains a commutative subalgebra ${\mathfrak h}$ of dimension from the third table row.
  
Since, by assumption, $\dim {\mathfrak h}\leq n$, it is clear that the inequality~(\ref{E:A}) holds for algebras of the first five types. It remains to check that it is true for infinite series:

Case 1: ${\mathfrak g}$ is of type $A_l,\ l\geq1,\ \dim {\mathfrak h}\geq 1$.
   $$  \dim {\mathfrak g} = (l+1)^2 -1 \leq  4\left[\frac{(l+1)^2}{4} \right]+3 -1= 4\dim {\mathfrak h} +2\leq 6\dim{\mathfrak h} < 7n .$$

Case 2: ${\mathfrak g}$  is of type $B_l,\ l\geq 3,\ \dim {\mathfrak h}\geq l$.
$$\dim {\mathfrak g} =2l^2+l =4\dim{\mathfrak h} +3l -4  < 7\dim{\mathfrak h}\leq 7n.$$

Case 3: ${\mathfrak g}$  is of type  $C_l, \ l\geq2$.
$$\dim {\mathfrak g} =2l^2+l \leq 2l(l+1)=4\dim{\mathfrak h} <7n.$$

Case 4: ${\mathfrak g}$ is of type $D_l,\ l\geq4,\ \dim {\mathfrak h}\geq l$.
  $$ \dim {\mathfrak g}=2l^2-l=4\dim{\mathfrak h}+l\leq 5\dim{\mathfrak h}<7n.\ \bl$$

\begin{lem} \label{T:C}
If ${\mathfrak g}$ is a complex semisimple Lie algebra satisfying the condition $\A$, then
$$\dim {\mathfrak g} \leq 7n.$$
\end{lem}
\doc
 Let ${\mathfrak g}={\mathfrak g}_1\oplus \dots \oplus {\mathfrak g}_m$ be the direct sum decomposition of  ${\mathfrak g}$ into simple subalgebras. For $ i=1,\dots, m$ we choose an abelian subalgebra ${\mathfrak h}_i \leq {\mathfrak g}_i$ of maximal dimension $n_i$. Then, by the previous Lemma, $\dim {\mathfrak g}_i \leq 7n_i$.
 
Further, ${\mathfrak h}={\mathfrak h}_1\oplus \dots \oplus {\mathfrak h}_m$ is abelian subalgebra of ${\mathfrak g}$. By assumption,  $\dim {\mathfrak h}=\sum_{i=1}^m n_i \leq n$. Thus we get  

$$ \dim {\mathfrak g} =  \sum_{i=1}^m \dim {\mathfrak g}_i \leq 7\sum_{i=1}^m n_i \leq 7n .\ \bl$$

\begin{thm} \label{T:F}
The function  $l_{{\C}}(n)$ satisfies the following inequality:
 $$ l_{{\C}}(n) \leq \frac{n^2+17n}{2} . $$

\end{thm}

\doc
Consider a complex Lie algebra ${\mathfrak g}$ satisfying the condition $\A$. Let us show that $\dim {\mathfrak g} \leq\frac{n^2+17n}{2}$.
 
By Levi's theorem, ${\mathfrak g}$ is the semidirect product
of a semisimple Lie algebra ${\mathfrak s}$ with the solvable radical ${\mathfrak r}$. 
Let  $n_1$ and  $n_2$ be the maximal dimensions of abelian subalgebras of ${\mathfrak r}$ and ${\mathfrak s}$ respectively.
Then, by Lemma~\ref{T:R}, $\dim {\mathfrak r} \leq \frac{n_1(n_1+3)} 2$ and, by Lemma~\ref{T:C}, $\dim {\mathfrak s} \leq 7n_2$.

  By assumption, $n_1,n_2\leq n$, so that 
$$ \dim {\mathfrak g} = \dim {\mathfrak r }+\dim {\mathfrak s } \leq \frac  {n_1(n_1+3)}{2} + 7n_2   \leq     \frac{n^2+17n}{2}.\      \bl$$

\begin{sled} \label{T:T}
For the function  $g_{\C}(n)$ the following inequality hols:
$$g_{\C}(n) \leq  \frac{n^2+17n}{2}.$$
\end{sled}
\doc
Let $G$ be a complex Lie group such  that the dimensions of all abelian Lie subgroups of $G$ are not greater than $n$, and let  ${\mathfrak g}$ be the Lie algebra associated with  $G$.
Consider an abelian subalgebra  ${\mathfrak h}\leq{\mathfrak g}$. We denote by ${\mathfrak h}^M$  the minimal subalgebra such that   ${\mathfrak h}\leq{\mathfrak h}^M$  and   there exists a connected Lie  subgroup $H\leq G$ with   Lie algebra ${\mathfrak h}^M$. By Theorem~1.4.3,~\cite{VO95}, the commutator subalgebras of ${\mathfrak h}$ and ${\mathfrak h}^M$ are equal. Therefore ${\mathfrak h}^M$ is commutative and hence so is $H$. By assumption, $\dim {\mathfrak h}\leq \dim {\mathfrak h}^M=\dim H\leq n$. We see that ${\mathfrak g}$ satisfies the condition $\A$.
Now, by Theorem~\ref{T:F}, we get
$\dim G=\dim {\mathfrak g} \leq \frac{n^2+17n}{2}.$~$\bl$

\subsection{Real Lie algebras}
\begin{lem} \label{T:S}
The dimension of a real solvable Lie algebra $\mathfrak g$ that satisfies the  condition $\A$ is not greater than $2n^2+3n.$
\end{lem}
\doc
Let $\mathfrak i\lhd \mathfrak g$ be a maximal ideal such that $\mathfrak i$ is commutative or $\mathfrak i$ is a sum of non-commutative ideals $\mathfrak i$$_k$ such that
\begin{enumerate}
\item [(i)] $ [{\mathfrak i}_k, {\mathfrak i}_l]=0,\ \mbox{for } k\ne l,\ k,l=1,\dots,m;$
\item  [(ii)]${\mathfrak i}_k \cap {\mathfrak i}_l={\mathfrak c},\ \mbox{for } k\ne l,\ k,l=1,\dots,m, $  where ${\mathfrak c}$ is a central ideal of ${\mathfrak i}$;
\item [(iii)]$\dim {\mathfrak i}_k /{\mathfrak c}=2,\ \mbox{for}\  k=1,\dots,m$.
\end{enumerate}
If ${\mathfrak i}$ is commutative put ${\mathfrak c}={\mathfrak i}$ and $m=0$.
   
For each $x\in{\mathfrak g}$ put $\rho(x)=\ad x|_{\mathfrak i}$. Let us show that $\Ker \rho ={\mathfrak c}$. 
 
Since all  ideals ${\mathfrak i}_k$ are non-commutative, it follows from (iii) that ${\mathfrak c}$  is a centre of each~${\mathfrak i}_k$.  Hence $\Ker \rho \cap {\mathfrak i}={\mathfrak c}$.  Suppose that ${\mathfrak c} \subsetneq \Ker \rho$. Then   $\Ker \rho / {\mathfrak c}$  contains a one- or two-dimensional ideal ${\mathfrak j} /{\mathfrak c}$ of the real solvable factor-algebra ${\mathfrak g}/{\mathfrak c}$ (see Corollary~4 of Theorem~I.5.1,~\cite{B76}). If ${\mathfrak j}$ is commutative put   ${\mathfrak c}'={\mathfrak j}, \ {\mathfrak i}_k'={\mathfrak i}_k+{\mathfrak j},\ k=1,\dots,m,$ and $ {\mathfrak i}'={\mathfrak i}'_1+\dots+{\mathfrak i}'_m$. It is clear that ${\mathfrak i}'$ is greater than  ${\mathfrak i}$ and that ${\mathfrak i}'$ satisfies all the conditions we need. If ${\mathfrak j}$ is not commutative, then, since ${\mathfrak c}$ belongs to the centre of ${\mathfrak j}$, $\dim {\mathfrak j}/{\mathfrak c}=2$. Putting ${\mathfrak i}_{m+1}={\mathfrak j}$ and $\ {\mathfrak i}'={\mathfrak i}+{\mathfrak i}_{m+1}$, we get the ideal ${\mathfrak i}'$ such that ${\mathfrak i} \subsetneq{\mathfrak i}'$ and ${\mathfrak i}'$ satisfies the same conditions  again. We have a contradiction.

Thus  ${\mathfrak g}/{\mathfrak c}$ is naturally embedded as a solvable subalgebra into ${\mathfrak gl}_s(\R)$ where $s=\dim {\mathfrak i}$.  Since the dimension of any irreducible representation of a solvable real Lie algebra is at most 2, we can choose a basis of ${\mathfrak i}$ such that all elements of ${\mathfrak g}/{\mathfrak c}$ are represented in this basis by "almost triangular"  matrices, that is these matrices have arbitrary $1\times1$ or $2\times2 $ submatrices   on the main diagonal and zeros only below the diagonal. It follows from this that $\dim {\mathfrak g}/{\mathfrak c}\leq \frac {s(s+2)} 2$.

Now notice that ${\mathfrak i}$  contains a "large" commutative  subalgebra.
Indeed,  by condition~(i),  for any  $k$  the set ${\mathfrak i}_k \cap \oplus\,_{l\ne k}{\mathfrak i}_l$ belongs to the centre of
${\mathfrak i}_k$ and hence this set coincides with ${\mathfrak c}$.  Consequently, ${\mathfrak i} / {\mathfrak c}$ is a direct sum of ${\mathfrak i}_1/ {\mathfrak c}, \dots, {\mathfrak i}_m / {\mathfrak c} $  and $\dim {\mathfrak i}=2m+\dim {\mathfrak c}$.
For $k=1,\dots,m$ we choose $x_k\in {\mathfrak i}_k\setminus {\mathfrak c}$. 
Then the subspace spanned by $  \itd{x}{m}$ and ${\mathfrak c} $ is abelian subalgebra of dimension $t=m+\dim{\mathfrak c}\geq \frac s 2$.
  We get
$$\dim {\mathfrak g} \leq \frac{s(s+2)} 2 + \dim {\mathfrak c}\leq t(2t+2)+t=2t^2+3t.\ \bl$$

\begin{lem} \label{T:G}
If ${\mathfrak g}$  is a real simple Lie algebra satisfying the  condition $\A$, then
 \begin{equation*}  
 \dim {\mathfrak g} \leq 2n^2 +15n.
\end{equation*}
\end{lem}
\doc
We denote by ${\mathfrak u}^{\R}$ the realification of a complex Lie algebra ${\mathfrak u}$.

The Theorem~5.1.1,~\cite{VO95} says that a real Lie algebra is simple if and only if either it is isomorphic to ${\mathfrak u}^{\R}$ for some complex simple Lie algebra $\mathfrak u$ or it is a real form of a complex simple Lie algebra.

First we assume that ${\mathfrak g}={\mathfrak u}^{\R}$ for some complex simple Lie algebra $\mathfrak u$. It follows from  Lemma~\ref{T:B} that  there exists an abelian subalgebra ${\mathfrak h}\leq{\mathfrak u}$ such that ${\dim_{\C} {\mathfrak u} \leq 7\dim_{\C} {\mathfrak h}}$. Then ${\mathfrak h}^{\R}$ is abelian subalgebra of ${\mathfrak u}^{\R}$ and hence, since ${\dim_{\R} {\mathfrak g}=2 \dim_{\C} {\mathfrak u}}$ è ${\dim_{\R} {\mathfrak h}^{\R} =2 \dim_{\C} {\mathfrak h}}$, in this case the Lemma is true.

Suppose now that ${\mathfrak g}$ is a real form of a complex simple algebra ${\mathfrak u}$. According to the Exercise~5.4.9,~\cite{VO95}, there exists a commutative subalgebra ${\mathfrak h}\leq{\mathfrak g}$ such that the complexification ${\mathfrak h}(\C)$ is  a Cartan subalgebra of ${\mathfrak u}$. In particular, the real dimension of ${\mathfrak h}$ equals the rank of~${\mathfrak u}$. Now, comparing the ranks and the dimensions of the simple complex Lie algebras from the Table, we can obtain easily the inequality we need.~$\bl$

\begin{lem} \label{T:H}
The dimension of a real semisimple Lie algebra ${\mathfrak g}$ that satisfies  the condition  $\A$ is majorized by $ 2n^2+ 15n. $
\end{lem}
\emph{The proof} is analogous to the proof of Lemma~\ref{T:C}. It is enough to note that the sum of  positive numbers squared is majorized by the square of the sum   of the numbers.~$\bl$

\begin{thm} \label{T:I}
For the function $l_{{\R}}(n)$ the following relation holds:
$$ l_{{\R}}(n) \leq 4n^2+18n.$$
\end{thm}
\emph{The proof} is similar to the proof of  Theorem~\ref{T:F}. The only difference is in that we use Lemmas~\ref{T:R} and~\ref{T:C} instead of Lemmas~\ref{T:S} and~\ref{T:H}.~$\bl$
 
\begin{sled}
The function $g_{\R}(n)$ satisfies the inequality
$$g_{\R}(n) \leq 4n^2+18n. $$
\end{sled}
\emph{The proof} is similar to the proof of Corollary~\ref{T:T}.~$\bl$

\subsection{Associative algebras}

\begin{lem} \label{T:N}
If $A$ is a nilpotent associative algebra over an arbitrary field such that $A$ satisfies the condition $\A$, then
$$\dim A \leq \frac{n(n+1)}{2}.$$
\end{lem}
\doc
We can consider $A$ as a nilpotent Lie algebra   with the commutator $[x,y]=xy-yx$ as product.
Let $B\leq A$ be a maximal abelian Lie subalgebra. For any  $x,y,z \in B$ we have
\begin{align*}
 xy=yx, \ \ &\mbox{since} \ [x,y]=0;\\  [xy,z]=0  , \ \ &\mbox{since} \  (xy)z= xzy=z(xy).
\end{align*}

Hence, since $B$ is maximal, $xy \in B$. Consequently $B$ is commutative associative subalgebra of $A$ and, by assumption, $\dim B \leq n$. Now we can use Lemma~\ref{T:A} to conclude the proof.~$\bl$
 
\begin{lem} \label{T:U}
The dimension of a semisimple associative algebra that satisfies the condition    $\A$ is not greater than $n^2$.
  \end{lem}
\doc
If $A$ is a simple algebra over a field $K$, then, by  Wedderburn's theorem, it is isomorphic to a matrix algebra $M_r(D)$ for some division algebra $D$ and positive integer~$r$.

Let $L$ be the centre of $D$, and let $F$ be a maximal subfield of $D$. By Theorem~4.2.2,~\cite{H72}, we have $[D:F]=[F:L].$  Therefore, observing $K\subset L $, we get
$$ [D:K]=[D:F][F:K]=[F:L][F:K] \leq [F:K]^2. $$
 
All diagonal matrices of $M_r(D)$ with entries in $F$ form a commutative subalgebra of dimension $s=r[F:K]$. It follows from this the estimate we need:
$$\dim A=r^2[D:K] \leq r^2 [F:K]^2 =s^2.$$

So, we proved the Lemma in the case of   a simple algebra $A$. Since any semisimple algebra is a direct sum of simple ones, the Lemma  is also true for any such an algebra.~$\bl$

\begin{lem} \label{T:J}
Suppose $K$ is a field satisfying one of the following conditions: 

 {\bfseries {\upshape (i)}}  the algebraic closure $\bar K$ is a finite extension of $K$ and
$c=\max\{[\bar K:K],\frac 9 2\}$; 

{\bfseries {\upshape (ii)}}  $K$ is finite and $c=\frac 9 2$.

 If a semisimple associative algebra over $K$ satisfies the condition $\A$, then its dimension   is majorized by $cn$.
  \end{lem}
\doc
First suppose that $A$ is a simple algebra. By Wedderburn's theorem, $A$ is isomorphic to a matrix algebra $M_r(D)$ for some division algebra $D$ and positive integer $r$.

As before, let $L$ be the centre of $D$, and let $F$ be a maximal subfield of $D$. By Theorem~4.2.2,~\cite{H72}, $[D:F]=[F:L].$
If $\bar K$ is a finite extension of $K$,  then, since $K\subset L \subset F \subset \bar K$, we have  $$[D:F]\leq c.$$
If $K$ is finite, then $D$ is also finite, so that ${[D:F]=1< c}$, since any finite division algebra is commutative (see Section~VI.3,~\cite{S83}).

Let $r=1$.  Then $A=D$ contains the commutative subalgebra $B=F$. We obtain  ${\dim A=[D:F][F:K] \leq c \dim B}$.

Further, let $r=2k$. Consider the abelian subalgebra $B\subset M_r(D)$, 
consisting of matrices whose only nonzero elements lie in a $k\times k$ submatrix in the upper right-hand corner.
 We have ${\dim B =k^2[D:K],}$ and  $$ \dim A =4k^2[D:K] \leq c \dim B. $$

Similarly, if $m=2k+1, \ k\geq 1$, the commutative subalgebra $B$  consisting of matrices whose only nonzero entries belongs to a $k\times (k+1)$ submatrix in the upper right-hand corner  has dimension ${k(k+1)[D:K]\geq 2[D:K]}$.
At the same time,
$$\dim A =(2k+1)^2[D:K]=(4k(k+1)+1)[D:K]\leq \frac{9\dim B} 2\leq c\dim B.$$
 
Since any semisimple algebra is a direct sum of simple ones, the general case follows from the previous arguments.~$\bl$

\begin {thm} \label{T:K}
Let $K$  be an algebraically closed field or a field of characteristic 0, then
$$a_K(n) \leq \frac{3n^2 +n}{2}.  $$
If, in addition, $c=max\{[\bar K:K],\frac 9 2\} <\infty$, then
$$a_{K}(n)\leq \frac {n^2+(2c+1)n} 2.$$
In particular, in the case of a real or   complex field we have
$$a_{{\R}}(n) ,a_{{\C}}(n)\leq \frac {n^2}{2}+5n.$$
\end{thm} 
\doc
Let $A$ be an algebra over $K$ and let $R$ be its Jacobson radical. The algebra $A/R$ is semisimple. If it is separable, the Wedderburn-Malcev theorem asserts that $A$ is isomorphic to the semidirect sum of $A/R$ and $R$. Note that since $K$ is algebraically closed or has    characteristic~0, $A/R$ is  separable.  Thus the first inequality follows from Lemmas~\ref{T:N} and~\ref{T:U} and the second one follows from Lemmas~\ref{T:N} and~\ref{T:J}.~$\bl$\\ \\
\textbf{Theorem~\ref{T:K}$'$.}
\emph{
Let $K$  be an algebraically closed field or a field of characteristic 0, then
$$a^1_K(n) \leq \frac{3n^2 +n}{2}.  $$
If, at the same time, $c=max\{[\bar K:K],\frac 9 2\} <\infty$, then }
$$a^{1}_{K}(n)\leq \frac {n^2+(2c+1)n} 2.$$
\doc 
Any maximal commutative subalgebra of a unital algebra is unital. Consequently if all commutative unital subalgebras of a unital algebra have dimensions at most $n$, then the algebra satisfies the condition $\A$. Hence this Theorem follows from the previous one.~$\bl$

\section{ Quadratic lower bounds}
\subsection{Lie algebras over an arbitrary field}
\begin{thm} \label{T:Q}
For any field  $K$ the following inequality holds:
$$ l_{K}(n) \geq \frac {n^2+4n-5}{8}.$$
\end{thm}

In order to prove the Theorem,  for each positive integer $s$ we shall construct a nilpotent Lie algebra of class 2 and of dimension at least $ \frac {s^2+4s-5}{8}$ such that any abelian subalgebra has  dimension at  most $s$. Here it is suitable to introduce special functions for nilpotent groups and algebras of class 2 that are analogous to the functions we defined before.

\begin{opr}
For any integer $n$ we denote by $l^n_{K}(n)$ ($a^n_{K}(n)$ and $g^{n}_K(n)$) the greatest integer $h$ such that there exists a nilpotent Lie algebra (an associative algebra and a Lie  group respectively) of class 2 and of dimension $h$ such that the dimensions of all commutative subalgebras (Lie subgroups) of this algebra (Lie group) are not greater than $n$.
\end{opr}

We shall prove the following statement, which is stronger than the Theorem~\ref{T:Q}.\\ \\
\textbf{Theorem~\ref{T:Q}$'$.}
\emph{
For any field $K$ we have the following estimate:
$$ l^{n}_{K}(n) \geq \frac {n^2+4n-5}{8}.$$
}

In the case of finite field $K$ the construction of examples of algebras is analogous to the construction used by Ol'shanskii~\cite{O78} for finite $p$-groups. And in the infinite case the way of construction is parallel to the way used in~\cite{M4} for  finitely generated torsion-free nilpotent groups.
\begin{lem} \label{T:Y}
 
    Suppose that the positive integers $k,\ t$, and $n$ satisfy  the inequality
\begin{equation*}  
  2n<t(k-1) ,
\end{equation*}
and let $V$ be an $n$-dimensional vector space over a field $K$. Then  there exists a ${t\mbox{-tuple}}$  \linebreak
 ${\varPhi=\{\ph_1, \dotsc, \ph_t \}}$  of skew-symmetric bilinear forms on   $V$   such that  no  $k$-dimensional subspace is simultaneously  isotropic for all of the forms   $\ph_1, \dotsc, \ph_t$. 
\end{lem}
\doc
This Lemma was proved in~\cite{O78} for a prime field $K$. But the same proof  obviously shows  that it is true if $K$ is an arbitrary finite filed. On the other hand, in the case of infinite field the statement of the Lemma is very close to   the Main Lemma of~\cite{M4}. (Here the inequality is stronger.) And  it is enough to modify the proof of that Lemma.

Indeed, let $\bar{K}$ be the algebraic closure  of $K$. Put $\bar{V} = V_K \otimes \bar K$.
Note that the   arguments used in the Main Lemma of~\cite{M4} hold true if the   field of complex numbers is replaced by an arbitrary closed field.  Arguing this way, it can be proved that the set of all $t$-tuples ${\varPhi=\{\ph_1, \dotsc, \ph_t \}}$ of skew-symmetric bilinear forms on $\bar V$ such that there exists a vector   subspace of dimension $k$ that is simultaneously  isotropic for all  $\ph_i$ is closed in Zariski topology. That is the set is given by a system of  equations that are polynomial in entries of matrices $\ph_1, \dotsc, \ph_t:$
$$G_i=0, \ \  i=1,\dots,r.$$
Moreover, there is $i$ such that $G_i \not\equiv 0.$
Therefore, since $K$ is infinite, the set does not contain all $t$-tuples $\varPhi$ with elements in $K$. This concludes the proof.~$\bl$

Let $U$ and $V$ be vector spaces over   $K$  having bases $\itd{f}{t}$ and $\itd{e}{n}$  respectively, and let 
$\varPhi=\{ \itd{\varphi}{t} \}$ be a $t$-tuple of skew-symmetric bilinear forms on $V$. Consider the direct sum $\gf$ of $U$ and $V$. Define the product of two elements  $x,y\in\gf$  by
\begin{equation} \label{E:h}
 [x,y]=\ph_1(v_1,v_2)f_1+\dots+ \ph_t(v_1,v_2)f_t    ,
\end{equation}
where $x=u_1+v_1,\ y=u_2+v_2$ for some $u_1,u_2\in U$, $v_1,v_2\in V$.
From~(\ref{E:h}) it follows that $\gf$ is a nilpotent Lie algebra of class 2 and that $U$ is a central $t$-dimensional subalgebra.
 
\begin{lem} \label{T:X}
        Suppose that the positive integers   $k,\ t$,  and $n$ satisfy the inequality   $2n< t(k-1)$,   and  that a $t$-tuple   $\varPhi=\{ \itd{\varphi}{t} \}$ of skew-symmetric  bilinear forms on an $n$-dimensional vector space $V$ is chosen in accordance with Lemma~\ref{T:Y}. Then the corresponding  Lie algebra  $ \gf$ has no  abelian subalgebras of dimension greater than  $ k+t-1$.
\end{lem}
\doc
Let $\mathfrak h$ be an abelian subalgebra of $\gf$. Then from~(\ref{E:h}) it follows that the subspace ${\mathfrak h/ (h} \cap U) \subset V$ is isotropic for $\itd{\ph}{t}$. Therefore, by Lemma~\ref{T:Y},
\[
   \dim  {\mathfrak h/ (h} \cap U) \leq   k-1.
\]
Thus,
   \[
 \dim  {\mathfrak h} \leq k-1+ \dim U = k+t-1. \ \bl
\]
 
\noindent \emph{Proof of Theorem~\ref{T:Q}$'$.}
Let $s$ be even. Put    $t=\frac{s}{2}+1$, $ k=\frac{s}{2}$,  $n=\left[ \frac{s^2-5}{8}\right]$. Then $2n< t(k-1)$ and, by Lemma~\ref{T:X}, there exists a Lie algebra   $\gf$ all of whose abelian subalgebras have dimensions at most  $s$ and such that the dimension of   $\gf$    is equal to

\begin{equation*} n+t=\left[\frac{s^2-5}{8}\right]+\frac{s}{2}+1 \geq\frac{s^2+4s-4}{8}. \end{equation*}

 \noindent Similarly, if $s$  is odd, putting $t=k=\frac{s+1}{2}$,   $n=\left[ \frac{s^2-2}{8}\right]$, we obtain that the dimension of  $\gf$ equals

\begin{equation*}
n+t=\left[\frac{s^2-2}{8}\right]+\frac{s+1}{2} \geq\frac{s^2+4s-5}{8}.
\end{equation*}

\noindent Hence $f(s)\geq\frac{s^2+4s-5}{8}$ and the Theorem is proved.~$\bl$

\subsection{Lower bounds for the functions $l_{\R}(n)$, $a_{K}(n)$, and $g_K(n)$ }
\begin{thm} \label{T:Z}
The function $l_{\R}(n)$ satisfies the relation  
\begin{equation} \label{E:Z}
l_{\R}(n)\geq 2n^2 +n.
\end{equation}
\end{thm}
\doc
For any positive integer $n$ consider  the real Lie algebra ${{\mathfrak g}_n={\mathfrak so}_{2n+1}(\R)}$ of real skew-symmetric matrices. According to the Example of the Subsection~4.1.3,~\cite{VO95},  it is compact real form of complex simple Lie algebra ${\mathfrak so}_{2n+1}(\C)$ of type $B_n$.
Let ${\mathfrak g}_n={\mathfrak f}\oplus{\mathfrak p}$ be a Cartan decomposition of  ${\mathfrak g}_n$ into the direct sum of the subalgebra ${\mathfrak f}$ and the subspace ${\mathfrak p}$. Since ${\mathfrak g}_n$ is compact, the decomposition is trivial, that is ${\mathfrak f}={\mathfrak g}_n$  and  ${\mathfrak p}=0$ (see Exercise~5.3.6,~\cite{VO95}). 
Further, consider a maximal commutative subalgebra ${\mathfrak h}$  of $ {\mathfrak f}$. Its complexification ${\mathfrak h}(\C)$ is a Cartan subalgebra of ${\mathfrak f}(\C)={\mathfrak g}$$_n(\C)={\mathfrak so}$$_{2n+1}(\C)$ (see Exercise~5.3.27,~\cite{VO95}). Thus we obtain that the dimension of any maximal abelian subalgebra of ${\mathfrak g}_n$ is equal to the rank of 
${\mathfrak so}_{2n+1}(\C)$, which equals  $n$.   At the same time, $\dim{\mathfrak g}_n=2n^2+n$. Therefore the inequality~\ref{E:Z} is proved.~$\bl$ 
 
\begin{note} \label{T:V}
We can get the same estimate if  instead of ${\mathfrak so}_{2n+1}(\R)$ we consider  the compact real form  ${\mathfrak sp}_n= \{\, \matr X Y {-\bar Y} {\bar X} \in {\mathfrak gl}_{2n}({\C}) \mid  {X,Y \in {\mathfrak gl}_n({\C}),} \ X=-\bar X ^{\top}, \ Y=Y^{\top} \}  $ of the complex simple Lie algebra $ {\mathfrak sp}_{2n}({\C})=  \{\, \matr X Y Z {-X^{\top}} \in {\mathfrak gl}_{2n}({\C}) \mid   $  $  { {X,Y,Z \in {\mathfrak gl}_n({\C}),} \ Y=-\bar Y ^{\top}, \ Z=Z^{\top} \}} $  of type $C_n$.
\end{note} 

\begin{thm} \label{T:O}
If a field $K$ has   characteristic other than 2, then  
$$ a^n_{K}(n) \geq \frac {n^2+4n-5}{8}.$$
\end{thm}
\doc
The proof of the Theorem~\ref{T:Q}$'$ shows that for any integer $n$ there exists a nilpotent  Lie algebra ${\mathfrak g}_n$  of class 2 and of dimension at least $ \frac {n^2+4n-5}{8}$ over $K$ such that any abelian subalgebra has  dimension at  most $n$. Since ${\mathfrak g}_n$ is nilpotent of class 2, the product of any three elements of ${\mathfrak g}_n$ equals 0. Hence ${\mathfrak g}_n$ is a nilpotent associative algebra of class two with respect to the same operations.
We note that, since the product in a Lie algebra is anti-symmetric and $K$ has characteristic other than 2, commutative associative subalgebras of ${\mathfrak g}_n$ are exactly abelian Lie subalgebras. Therefore the algebras  ${\mathfrak g}_n$  serve  as examples of associative class two nilpotent algebras of "large" dimensions with "small" dimensions of commutative subalgebras.~$\bl$ 

\begin{note} \label{T:W}
Theorem~\ref{T:O} can be proved without  using the result about Lie algebras, but applying similar arguments. Then the assumption that the field has   characteristic other than 2 is not necessary. Indeed, let $\varPhi=\{ \itd{\varphi}{t} \}$ be a $t$-tuple of arbitrary bilinear forms on a vector space $V$ of dimension $n$, and let $A_{\varPhi}$ be the direct sum of vector spaces $U$ and $V$.
Define the product of two elements of $A_{\varPhi}$ by 
 $$ (u_1+v_1)\centerdot(u_2+v_2)=\ph_1(v_1,v_2)f_1+\dots+ \ph_t(v_1,v_2)f_t   , $$
where $u_1,u_2 \in U,\ v_1,v_2\in V$, and $\itd{f}{t}$ is a basis of $U$.  Then $A_{\varPhi}$ is a nilpotent associative algebra of class 2. And it follows from this formula  that a subalgebra $B$ of $A_{\varPhi}$ is commutative if and only if 
the restrictions of    $\itd{\varphi}{t}$  to $B /(B\cap U) \subset V$  are symmetric.
Further, arguing as in Lemma~\ref{T:Y},  it can be shown that if the integers $k,\ t$  and  $n$ satisfies the same inequality, then there exists a $t$-tuple $\varPhi=\{ \itd{\varphi}{t} \}$ of bilinear  forms on $V$ such that their restrictions to any $k$-dimensional subspace are not simultaneously symmetric. The corresponding algebra $A_{\varPhi}$ for this $t$-tuple has no commutative subalgebras of dimension greater than $k+t-1$. Finally, putting $k,\ t$ and  $n$ as in the proof of Theorem~\ref{T:Q}$'$, we get the same estimate.
\end{note}

  Since, by definition, $a_K(n)\geq a_K^n(n)$, the lower bound for  the  function $a_K(n)$ from the Main Theorem follows from this Remark. \\ \\
\textbf{Theorem~\ref{T:O}$'$.}
\emph{For any field $K$ the following holds:
\begin{equation} \label{E:Y}
 a^{1}_{K}(n) \geq \frac {n^2+2n}{8}.
\end{equation}
}
\doc
 Consider an associative algebra $A$.  $A$ is embedded in a unital associative algebra $A_1$ as follows. Let $A_1$ be the direct sum of $A$ and  a one-dimensional vector space spanned by  a vector $e$, i.e. $A_1=A\oplus Ke$.
Extend the operation of product from $A$ to $A_1$ according to the formula $xe=ex=x$ for  $x\in A_1$ and by linearity.
By definition, $e$ is an identity element of~$A_1$. Any maximal commutative subalgebra of $A_1$ has the form $B_1=B\oplus Ke$ for some maximal commutative subalgebra $B \leq A$. Hence we have $\dim B_1=\dim B+1$ and $\dim A_1=\dim A+1$. 

In other words, if there exists an associative algebra of dimension $h$ satisfying the condition ${{\mathcal{A}}(n-1)}$, then there is  a unital associative algebra of dimension $h+1$ whose unital commutative subalgebras all have dimension at most $n$.
In our notations,
$$ a^{1}_{K}(n)\geq a_{K}(n-1)+1\geq a^n_{K}(n-1)+1, \ \ \mbox{for } n\geq 2. $$
Using the previous Theorem and   Remark~\ref{T:W}, we get (\ref{E:Y}) for  $n\geq 2$. Clearly, (\ref{E:Y}) is also true for $n=1$.~$\bl$
\begin{thm} \label{T:P}
 If $K$ is a real or  complex number field, then
$$ g_{K}(n) \geq g_K^n(n)\geq \frac {n^2+4n-5}{8}.$$
Also, we have
$$g_{\R}(n) \geq 2n^2 +n.$$
\end{thm}
\doc
The Theorem follows easily from Theorems~\ref{T:Q} and~\ref{T:Z} and from the following facts: given a complex or real Lie algebra, there exists a connected Lie group with this Lie algebra (see Theorem~6.2,~\cite{VO95}); a connected Lie group with a nilpotent Lie algebra is nilpotent of the same class; and commutative Lie subgroups correspond to commutative Lie subalgebras.~$\bl$

\section{On the nilpotent groups and algebras of class 2}
Lemmas~\ref{T:A} and~\ref{T:N} show that the dimension  of a nilpotent algebra that satisfies the condition $\A$ is asymptotically bounded above by $\frac {n^2} 2$. In fact, it  turns out that the restrictions  of the functions we study to nilpotent algebras and groups of class two still have quadratic growth.   Moreover, upper bounds  are asymptotically two times better, so that our estimates for the functions  we defined in  section~2 are quite precise.

\begin{thm}
If   $f$ is one of the functions $l^n_K,\ a^n_K$ or $g_K^n$, then it  satisfies the inequality
 $$\frac {n^2+4n-5}{8} \leq f(n)\leq \frac{n^2}4+n.$$
\end{thm}
\doc
The first inequality for the functions $l^n_K,\ a^n_K$ and $g_K^n$ is proved in  Theorems~\ref{T:Q}$'$,~\ref{T:O}  and~\ref{T:P} respectively. Let us show the second one.
 
Consider a nilpotent Lie algebra ${\mathfrak g}$ of class 2. Let ${\mathfrak z}$ be the centre of  ${\mathfrak g}$,  $\itd{z}{t}$ a basis of ${\mathfrak z}$, and $V$ a complementary subspace to ${\mathfrak z}$ of dimension $m$. Then the product of two elements $x=\bar{x}+\bar{\bar{x}}$ and $y=\bar{y}+\bar{\bar{y}}$ of ${\mathfrak g}$ with $\bar{x},\bar{y}\in V$ and  $\bar{\bar{x}},\bar{\bar{y}}\in{\mathfrak z}$ has the form 
 $$ [x,y]=\ph_1(\bar{x},\bar{y})z_1+ \dots+ \ph_t(\bar{x},\bar{y})z_t    $$
for some skew-symmetric bilinear forms  $\itd{\ph}{t}$ on $V$.

In this notation, $x$ and $y$ commute if and only if 
 $$ \ph_i(\bar{x},\bar{y})=0,\ \ \ \mbox{for } i=1,\dots, t. $$

Hence, given a nonzero element  $\bar{x}_1\in V$, the set of all elements $\bar{y}\in V$ that commute with $\bar{x}_1$ is given by the system of $t$ linear equations. Whence this set is a vector space of dimension at least $m-t$.  Consequently if ${m-t\geq2}$, then there is $\bar{x}_2 \in V$ such that $\bar{x}_1$ and $\bar{x}_2$ are linearly independent and $[\bar{x}_1,\bar{x}_2]=0$.

Further, proceeding by induction, suppose one has already chosen $k-1$ linearly independent elements $\itd{\bar{x}}{k-1}\in V$ such that $[\bar{x}_i,\bar{x}_j]=0$ for $i,j=1,\dots,k-1$. Then the set of all elements of $V$ that commute with any $\bar{x}_i$ is given by the system of $(k-1)t$ linear equations
$$ \ph_i(\bar{x_j},\bar{y})=0,\ \ \ \mbox{for }i=1,\dots, t,\  j=1,\dots,k-1. $$
Hence this set has dimension at least $m-(k-1)t$. Now if $m-(k-1)t\geq k$, then there exist $k$ linearly independent elements $\itd{\bar{x}}{k}\in V$ that commute with each other. Thus ${\mathfrak g}$ has an abelian subalgebra of dimension  $s=k+t$, which is generated by $\itd{\bar{x}}{k}$ and~${\mathfrak z}$.

Choose maximal $k$ satisfying this condition, that is  such that $m-(k-1)t\geq k$ and $m-kt\leq  k.$ Then 
$$ \dim  {\mathfrak g}=m+t\leq kt+k+t=(s-t)t+s \leq \frac{s^2}4+s .   $$

 Thus we proved the Theorem for $f=l^n_K$. If $f=a^n_K$, the Theorem follows from the previous arguments as Lemma~\ref{T:N} follows from Lemma~\ref{T:A}. And if $f=g^n_K$, using the previous arguments and taking into account that a nilpotent Lie group of class 2 has a class two nilpotent Lie algebra, we get the Theorem just as we proved Corollary~\ref{T:T}  using Theorem~\ref{T:F}.~$\bl$\\
 
\textbf{Acknowledgment. } The author would like to express her gratitude to Professor A.Yu.~Olshanskii  for  suggesting the problem,  his constant attention to this work, and his valuable advice and to A.~Minasyan and I.V.~Arzhantsev for  helpful comments.

{\scriptsize Higher Algebra, Department of Mechanics nd Mathematics, M.V.~Lomonosov Moscow State University, Leninskie Gory, GSP-2, Moscow; 119992. }

{\scriptsize  E-mail: mariamil@yandex.ru}

\end{document}